\documentclass[a4paper,12pt]{amsart}

\usepackage{cite,amsmath,amssymb,amsthm}
\usepackage[top=2.5cm,bottom=2.5cm,left=3.0cm,right=3.0cm]{geometry} 

\usepackage[margin=1cm,%
           font=small,%
           format=hang,%
           labelsep=period,%
           labelfont=bf]{caption}
\usepackage{graphicx}

\usepackage{hyperref}
\hypersetup{
	colorlinks   = true, 
	urlcolor     = blue, 
	linkcolor    = blue, 
	citecolor   = blue 
}

\usepackage{pgf,tikz,pgfplots}
\usepackage{mathrsfs}

\usepackage{amsmath,amsfonts,amssymb,amsthm,mathtools} 

\newtheorem{proposition}{Proposition}[]
\newtheorem{problem}{Problem}[]

\renewcommand{\baselinestretch}{1.2}

\begin{document}

\title{On the Wiener (r,s)-complexity of fullerene graphs}

\author{Andrey~A.~Dobrynin} 
\address{Sobolev Institute of Mathematics, Siberian Branch of the
Russian Academy of Sciences, Novosibirsk, 630090, Russia} 
\email{dobr@math.nsc.ru}

\author{Andrei~Yu.~Vesnin}
\address{Sobolev Institute of Mathematics, Siberian Branch of the
Russian Academy of Sciences, Novosibirsk, 630090, Russia   \newline \indent 
Tomsk State University, Tomsk, 634050, Russia}
\email{vesnin@math.nsc.ru}

\bigskip

\begin{abstract}
Fullerene graphs are mathematical models of fullerene molecules.
The Wiener $(r,s)$-complexity of a fullerene graph $G$ with vertex set $V(G)$
 is the number of pairwise distinct values of
$(r,s)$-transmission $tr_{r,s}(v)$ of its vertices $v$: 
$tr_{r,s}(v)= \sum_{u \in V(G)} \sum_{i=r}^{s} d(v,u)^i$
for positive integer $r$ and $s$.
The Wiener $(1,1)$-complexity is known as the Wiener complexity of a graph.
Irregular graphs have maximum complexity equal to the number of vertices.
No irregular fullerene graphs are known for the Wiener complexity.
Fullerene (IPR fullerene) graphs with $n$ vertices
having the maximal Wiener $(r,s)$-complexity are counted
for all $n\le 100$ ($n\le 136$) and small $r$ and $s$.
The irregular fullerene graphs are also presented.
\end{abstract}




\maketitle

\section*{Introduction}

A fullerene is a spherically shaped molecule consisting of carbon atoms
in which every carbon ring forms a pentagon or a hexagon
\cite{Fowl95,Schw15}.
Every atom of a fullerene has  bonds  with  exactly  three  neighboring atoms.
Fullerenes are the subject of intense research in chemistry and
they have found promising technological applications,
especially in nanotechnology and materials science
\cite{Ashr16,Cata11}.
Molecular graphs of fullerenes are called
\emph{fullerene graphs}.
A fullerene graph is a \mbox{3-connected} planar graph in which every vertex has degree 3
and every face has size 5 or 6.
By Euler's polyhedral formula, the number of pentagonal faces is always 12.
It is known that fullerene graphs having $n$ vertices exist for all even $n \ge 24$ and for $n = 20$.
The number of all non-isomorphic fullerenes was reported, for example, in
\cite{Brin97,Fowl95,Goed15-2}.
Fullerenes without adjacent pentagons,
\emph{i.\,e.},
each pentagon is surrounded
only by hexagons, satisfy the isolated pentagon rule and
called
\emph{IPR fullerenes}.
They are considered as  thermodynamic stable fullerene compounds.
Mathematical studies of fullerenes include
applications of graph theory and topology methods,
design of computational and combinatorial algorithms,
information theory approaches,
etc.
(see selected publications
\cite{Aliz14,Ando16,Ashr16,Cata11,Dobr20,Egor20,Fowl95,Ghor20,Goed15-2,Sabi18,Sabi15,Schw15}).

Topological indices are often applied for
quantifying structural complexity of molecular graphs.
One of the recent concept of complexity
reflects the diversity of values of a vertex graph invariant.
A number of distance topological indices of a graph are designed as
the sum of contributions of its vertices.
\emph{Complexity}
of a graph is
the number of pairwise distinct values of the vertex contributions.
Graphs with maximum possible complexity (equal to the number of vertices)
is called
\emph{irregular graphs}.
An irregular graph has obviously the identity automorphism group
while a graph with such a group may not be irregular.

One of the famous distance topological index is the Wiener index
of a graph $G$ with vertex set $V(G)$ which is defined as
$W(G)= \frac{1}{2}\sum_{v \in V(G)} tr(v)$,
where
$tr(v)= \sum_{u \in V(G)} d(v,u)$ is the transmission
of vertex $v$ and
$d(v,u)$ denotes the shortest distance between vertices $v$ and $u$.
This index was introduced by Harold Wiener in 1947 \cite{Wien47}.
Bibliography on the Wiener index and its applications can be found in
\cite{Dobr01,Dobr02,Knor16}.
The Wiener complexity of a graph is the number of pairwise distinct
vertex transmissions.
It was studied for various classes of graphs in
\cite{Aliz14,Aliz18,Dobr19-3,Dobr19}.
In particular, it was demonstrated that there do not exist
transmission irregular fullerene graphs with $n \le 232$ vertices
\cite{Dobr19}.
One of the possible reason of this fact is that
the interval of transmission values
may be too narrow for fullerene of such sizes.
In this paper, we consider a generalization of the vertex transmission
and compute the corresponding complexity of fullerene graphs.
Our goal is to find irregular fullerene graph.

\section*{Computational results}

Values of vertex transmission of arbitrary connected graphs with $n$ vertices
lie between $n-1$ and $n(n-1)/2$.
A generalization of the transmission $tr(v)$ of a vertex $v\in V(G)$ is the $(r,s)$-transmission
of $v$ defined as
$$
tr_{r,s}(v)= \sum_{u \in V(G)} \sum_{i=r}^{s} d(v,u)^i.
$$
For $r=s=1$, transmission $tr_{r,s}(v)$ coincides with $tr(v)$.
The interval of possible values of $tr_{r,s}(v)$
is much lager than that for $tr(v)$.
The corresponding analogs of the Wiener index of a graph $G$ may be written as
$$
W_{r,s}(G)\ = \sum_{\{u,v\} \subseteq V(G)} \sum_{i=r}^{s} d(u,v)^i = \frac{1}{2} \sum_{v \in V(G)} tr_{r,s}(v).
$$
The similar modifications of the Wiener index were considered in
\cite{Gutm97,Gutm04,Shab13}.
It can be noted that
$W_{1,1}$ is the Wiener index,
$W_{k,k}$ is the $k$-th distance moment for positive integer $k$
\cite{Klei99},
$WW = \frac{1}{2} W_{1,2} $ is the hyper-Wiener index
\cite{Rand93,Klei95},
and
$TSZ = \frac{1}{6}(2W_{1,2} + W_{2,3})$ is the Tratch--Stankevich--Zefirov index
\cite{Trat90}.

We have computed Wiener $(r,s)$-complexity for small values of $r$ and $s$, namely,
for $r,s=1,2,3$.
One of our goal is to find  irregular  fullerene graphs.
Computing results for  fullerene and  IPR fullerene graphs are presented in
Table~\ref{Table_1} and Table~\ref{Table_2},
respectively.
Here $C_{n,r,s}$ denotes the maximal
Wiener $(r,s)$-complexity
 of fullerene graphs with $n$ vertices.
The numbers of fullerene graphs with $C_{n,r,s}=n$,
i.e. with the maximal Wiener $(r,s)$-complexity, are in bold.
 Tables~\ref{Table_1} and \ref{Table_2}
show that irregular fullerene
graphs with $n\ge 76$ vertices exist for all considered $(r,s)$-transmissions.

{\renewcommand{\baselinestretch}{1}
{
\begin{table}[h!]
\centering
\caption{The number of fullerene graphs ($N$) with $n$ vertices having the
         maximal Wiener $(r,s)$-complexity $C_{n,r,s}$.} \label{Table_1}
\begin{tabular}{rr@{\hspace{1mm}}rr@{\hspace{1mm}}rr@{\hspace{1mm}}rr@{\hspace{1mm}}rr@{\hspace{1mm}}rr@{\hspace{1mm}}r} \hline
    & \multicolumn{2}{c}{$r=s=1$}
    & \multicolumn{2}{c}{$r=s=2$}
    & \multicolumn{2}{c}{$r=1, s=2$}
    &  \multicolumn{2}{c}{$r=s=3$}
    &  \multicolumn{2}{c}{$r=2, s=3$}
    &  \multicolumn{2}{c}{$r=1, s=3$}  \\ \hline
$n$ &$C_{n,r,s}$& $N$ &$C_{n,r,s}$ & $N$  & $C_{n,r,s}$ & $N$ & $C_{n,r,s}$& $N$ & $C_{n,r,s}$ & $N$ & $C_{n,r,s}$ & $N$ \\ \hline
 20 &   1 &  1  &    1  & 1    &    1    & 1 &     1  &  1     &   1  & 1       &  1      & 1       \\
 24 &   2 &  1  &    2  & 1    &    2    & 1 &     2  &  1     &   2  & 1       &  2      & 1       \\
 26 &   2 &  1  &    3  & 1    &    3    & 1 &     3  &  1     &   3  & 1       &  3      & 1       \\
 28 &   5 &  1  &    5  & 1    &    5    & 1 &     5  &  1     &   5  & 1       &  5      & 1       \\
 30 &   7 &  1  &    8  & 1    &    8    & 1 &     8  &  1     &   8  & 1       &  8      & 1       \\
 32 &   9 &  1  &    11 & 1    &   11    & 1 &     11 &  1     &  11  & 1       & 11      & 1       \\
 34 &  10 &  2  &    12 & 1    &   12    & 1 &     12 &  1     &  12  & 1       & 12      & 1       \\
 36 &  14 &  1  &    20 & 1    &   20    & 1 &     20 &  1     &  20  & 1       & 20      & 1       \\
 38 &  18 &  1  &    22 & 1    &   22    & 1 &     22 &  1     &  22  & 1       & 22      & 1       \\
 40 &  19 &  1  &    25 & 1    &   25    & 1 &     25 &  1     &  25  & 1       & 25      & 1       \\
 42 &  22 &  1  &    30 & 1    &   30    & 1 &     30 &  1     &  30  & 1       & 30      & 1       \\
 44 &  25 &  1  &    32 & 1    &   32    & 1 &     32 &  2     &  32  & 2       & 32      & 2       \\
 46 &  25 &  4  &    37 & 1    &   37    & 1 &     37 &  1     &  37  & 1       & 37      & 1       \\
 48 &  30 &  1  &    41 & 1    &   39    & 1 &     41 &  1     &  41  & 1       & 41      & 1       \\
 50 &  35 &  1  &    42 & 2    &   42    & 2 &     42 &  4     &  42  & 2       & 42      & 4       \\
 52 &  36 &  1  &    46 & 2    &   46    & 1 &     47 &  1     &  46  & 2       & 47      & 1       \\
 54 &  37 &  1  &    50 & 2    &   48    & 2 &     51 &  1     &  51  & 1       & 51      & 1       \\
 56 &  40 &  1  &    52 & 1    &   51    & 2 &     52 &  3     &  52  & 3       & 52      & 3       \\
 58 &  43 &  2  &    57 & 1    &   55    & 1 &     57 &  1     &  56  & 2       & 57      & 1       \\
 60 &  44 &  3  &    58 & 1    &   57    & 2 &     58 &  3     &  58  & 2       & 58      & 3       \\     
\bf{62} &  46 &  3  &    60 & 2    &   59    & 4 & \bf{62}&  1     &\bf{62} & 1     &\bf{62}  & 1       \\ 
\bf{64} &  49 &  5  &\bf{64}& 1    &   63    & 1 & \bf{64}&  2     &\bf{64} & 2     &\bf{64}  & 2       \\ 
\bf{66} &  50 &  2  &    65 & 5    &   65    & 1 & \bf{66}&  2     &\bf{66} & 2     &\bf{66}  & 2       \\
\bf{68} &  56 &  1  &    67 & 5    &   66    & 7 & \bf{68}&  1     &\bf{68} & 1     &\bf{68}  & 1       \\
\bf{70} &  56 &  1  &    69 & 7    &   69    & 1 & \bf{70}&  9     &\bf{70} & 9     &\bf{70}  & 9       \\ 
\bf{72} &  56 &  6  &\bf{72}& 2    &   71    & 4 & \bf{72}&  18    &\bf{72} & 16    &\bf{72}  & 18      \\
\bf{74} &  61 &  1  &\bf{74}& 4    &   73    & 3 & \bf{74}&  24    &\bf{74} & 19    &\bf{74}  & 26      \\ 
\bf{76} &  63 &  1  &\bf{76}& 2    &\bf{76}  & 1 & \bf{76}&  53    &\bf{76} & 50    &\bf{76}  & 55      \\
\bf{78} &  64 &  2  &\bf{78}& 14   &\bf{78}  & 1 & \bf{78}&  86    &\bf{78} & 72    &\bf{78}  & 92      \\
\bf{80} &  66 &  2  &\bf{80}& 14   &\bf{80}  & 2 & \bf{80}&  169   &\bf{80} & 140   &\bf{80}  & 174     \\
\bf{82} &  71 &  1  &\bf{82}& 22   &\bf{82}  & 3 & \bf{82}&  286   &\bf{82} & 251   &\bf{82}  & 299     \\
\bf{84} &  70 &  2  &\bf{84}& 52   &\bf{84}  & 11& \bf{84}&  483   &\bf{84} & 416   &\bf{84}  & 505     \\
\bf{86} &  73 &  3  &\bf{86}& 69   &\bf{86}  & 14& \bf{86}&  818   &\bf{86} & 672   &\bf{86}  & 856     \\
\bf{88} &  73 &  7  &\bf{88}& 132  &\bf{88}  & 16& \bf{88}&  1305  &\bf{88} & 1058  &\bf{88}  & 1345    \\
\bf{90} &  79 &  1  &\bf{90}& 154  &\bf{90}  & 36& \bf{90}&  2024  &\bf{90} & 1641  &\bf{90}  & 2104    \\
\bf{92} &  80 &  1  &\bf{92}& 247  &\bf{92}  & 38& \bf{92}&  3108  &\bf{92} & 2472  &\bf{92}  & 3292    \\
\bf{94} &  82 &  1  &\bf{94}& 385  &\bf{94}  & 73& \bf{94}&  4836  &\bf{94} & 3782  &\bf{94}  & 5052    \\
\bf{96} &  84 &  2  &\bf{96}& 511  &\bf{96}  & 86& \bf{96}&  6932  &\bf{96} & 5396  &\bf{96}  & 7366    \\
\bf{98} &  86 &  1  &\bf{98}& 697  &\bf{98}  &111& \bf{98}&  9800  &\bf{98} & 7623  &\bf{98}  & 10493   \\
\bf{100} &  89 &  1 &\bf{100}& 923  &\bf{100} &147&\bf{100}&  13870 &\bf{100}& 10627 &\bf{100} & 14886   \\  \hline
\end{tabular}
\end{table}
}

{\renewcommand{\baselinestretch}{1}
\begin{table}[h]
\centering
\caption{The number of IPR fullerene graphs ($N$) with $n$ vertices having
        the maximal Wiener $(r,s)$-complexity $C_{n,r,s}$.} \label{Table_2}
\begin{tabular}{rrrrrrrrrrrrr} \hline
    & \multicolumn{2}{c}{$r=s=1$}
    & \multicolumn{2}{c}{$r=s=2$}
    & \multicolumn{2}{c}{$r=1, s=2$}
    &  \multicolumn{2}{c}{$r=s=3$}
    &  \multicolumn{2}{c}{$r=2, s=3$}
    &  \multicolumn{2}{c}{$r=1, s=3$}  \\ \hline
$n$ &$C_{n,r,s}$& $N$ &$C_{n,r,s}$ & $N$ & $C_{n,r,s}$ & $N$ & $C_{n,r,s}$& $N$ & $C_{n,r,s}$ & $N$& $C_{n,r,s}$ & $N$ \\ \hline
 60 &  1  &  1  &    1   &  1   &      1    &  1  &    1   &   1    &     1   &    1   &     1     &   1     \\
 70 &  5  &  1  &    5   &  1   &     5     &  1  &    5   &   1    &     5    &   1   &     5     &   1     \\
 72 &  4  &  1  &    4   &  1   &     4     &  1  &    4   &   1    &     4    &   1   &     4     &   1     \\
 74 &  6  &  1  &    9   &  1   &     9     &  1  &    9   &   1    &     9    &   1   &     9     &   1     \\
 76 &  13 &  1  &    17  &  1   &     17    &  1  &    17  &   1    &     17   &   1   &     17    &   1     \\
 78 &  14 &  1  &    17  &  2   &     17    &  2  &    18  &   1    &     18   &   1   &     18    &   1     \\
 80 &  17 &  1  &    20  &  1   &     19    &  1  &    21  &   1    &     21   &   1   &     21    &   1     \\
 82 &  19 &  1  &    32  &  1   &     30    &  1  &    33  &   1    &     33   &   1   &     33    &   1     \\
 84 &  25 &  1  &    37  &  1   &     35    &  1  &    37  &   1    &     37   &   1   &     37    &   1     \\
 86 &  39 &  1  &    62  &  1   &     58    &  3  &    65  &   2    &     64   &   2   &     65    &   2     \\
 88 &  36 &  1  &    69  &  1   &     65    &  1  &    70  &   1    &     70   &   1   &     70    &   1     \\
 90 &  39 &  2  &    71  &  1   &     69    &  1  &    73  &   2    &     73   &   1   &     73    &   2     \\
 92 &  41 &  1  &    80  &  1   &     76    &  1  &    84  &   1    &     84   &   1   &     84    &   1     \\
 94 &  48 &  1  &    82  &  1   &     80    &  2  &    84  &   2    &     84   &   2   &     84    &   2     \\
 96 &  49 &  1  &    85  &  1   &     81    &  2  &    87  &   5    &     87   &   3   &     87    &   6     \\
 98 &  55 &  1  &    87  &  1   &     85    &  2  &    91  &   3    &     91   &   2   &     91    &   3     \\
100 &  55 &  3  &    90  &  2   &     90    &  1  &    95  &   1    &     95   &   1   &     95    &   1     \\
102 &  59 &  1  &    94  &  4   &     92    &  1  &    98  &   2    &     98   &   1   &     98    &   2     \\
104 &  65 &  1  &    97  &  2   &     99    &  1  &    101 &   1    &     100  &   9   &     101   &   2     \\
106 &  69 &  1  &    101 &  2   &     101   &  1  &    104 &   4    &     103  &   9   &     104   &   3     \\
108 &  70 &  1  &    103 &  2   &     102   &  3  &    107 &   3    &     107  &   2   &     107   &   2     \\ 
\bf{110} &  72 &  1  &    108 &  1   &     106   &  2  &\bf{110}&   1    & \bf{110} &   1   & \bf{110}  &   1     \\
\bf{112} &  74 &  1  &    108 &  3   &     109   &  2  &\bf{112}&   1    & \bf{111} &   8   & \bf{112}  &   1     \\
\bf{114} &  76 &  2  &    112 &  4   &     110   &  1  &\bf{114}&   4    & \bf{114} &   2   & \bf{114}  &   4     \\
\bf{116} &  80 &  2  &    114 &  2   &     113   &  1  &\bf{116}&  11    & \bf{116} &   6   & \bf{116}  &   15    \\
\bf{118} &  81 &  1  &    116 &  2   &     117   &  1  &\bf{118}&  32    & \bf{118} &   26  & \bf{118}  &   32    \\
\bf{120} &  87 &  1  &    119 &  1   &     118   &  2  &\bf{120}&  39    & \bf{120} &   34  & \bf{120}  &   47    \\
\bf{122} &  89 &  2  &    121 &  4   &     120   &  4  &\bf{122}&  73    & \bf{122} &   49  & \bf{122}  &   82    \\ 
\bf{124} &  91 &  2  &\bf{124}&  1   &     122   &  8  &\bf{124}&  146   & \bf{124} &   100 & \bf{124}  &   164   \\
\bf{126} & 93  &  1  &\bf{126}&  1   &     125   &  1  &\bf{126}&  262   & \bf{126} &   164 & \bf{126}  &   268   \\
\bf{128} & 96  &  1  &\bf{128}&  6   &     127   &  4  &\bf{128}&  409   & \bf{128} &   270 & \bf{128}  &   416   \\ 
\bf{130} &  97 &  3  &\bf{130}&  4   & \bf{130}  &  1  &\bf{130}&  739   & \bf{130} &   466 & \bf{130}  &   728   \\
\bf{132} & 100 &  2  &\bf{132}&  5   & \bf{132}  &  1  &\bf{132}&  1246  & \bf{132} &   749 & \bf{132}  &   1235  \\
\bf{134} & 102 &  2  &\bf{134}&  14  & \bf{134}  &  3  &\bf{134}&  2000  & \bf{134} &   1314& \bf{134}  &   1929  \\
\bf{136} & 105 &  1  &\bf{136}&  18  & \bf{136}  &  1  &\bf{136}&  3020  & \bf{136} &   1831& \bf{136}  &   2966  \\ \hline
\end{tabular}
\end{table}
}

\clearpage 

\begin{figure}[ht]
\begin{minipage}[h]{0.47\linewidth}
\center{\includegraphics[width=0.85\linewidth]{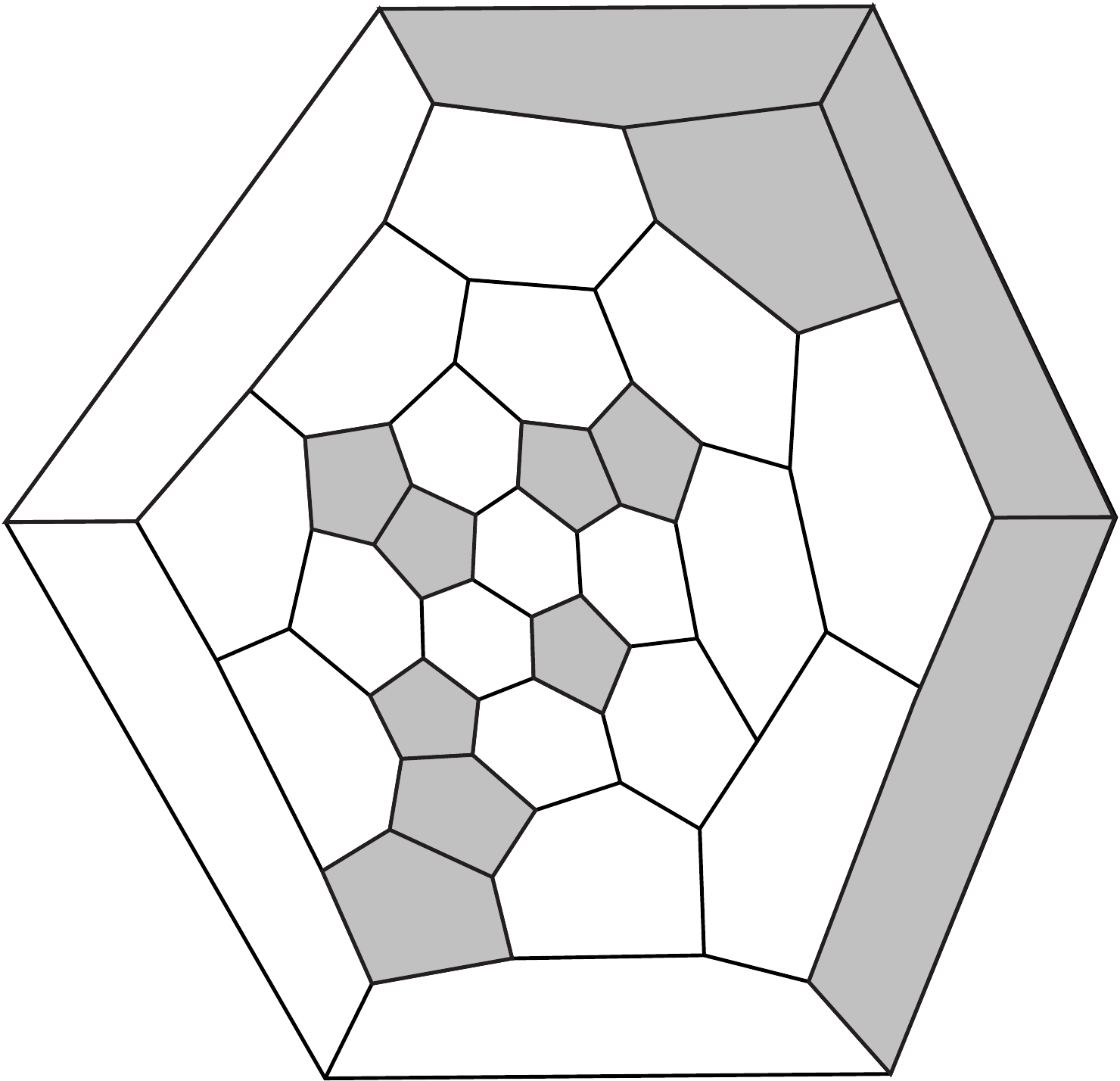}} $n=62 \ (r=1,2,3, s=3)$ \\
\vspace*{5mm}
\end{minipage}
\hfill
\begin{minipage}[h]{0.47\linewidth}
\vspace*{-5mm}
\center{\includegraphics[width=0.81\linewidth]{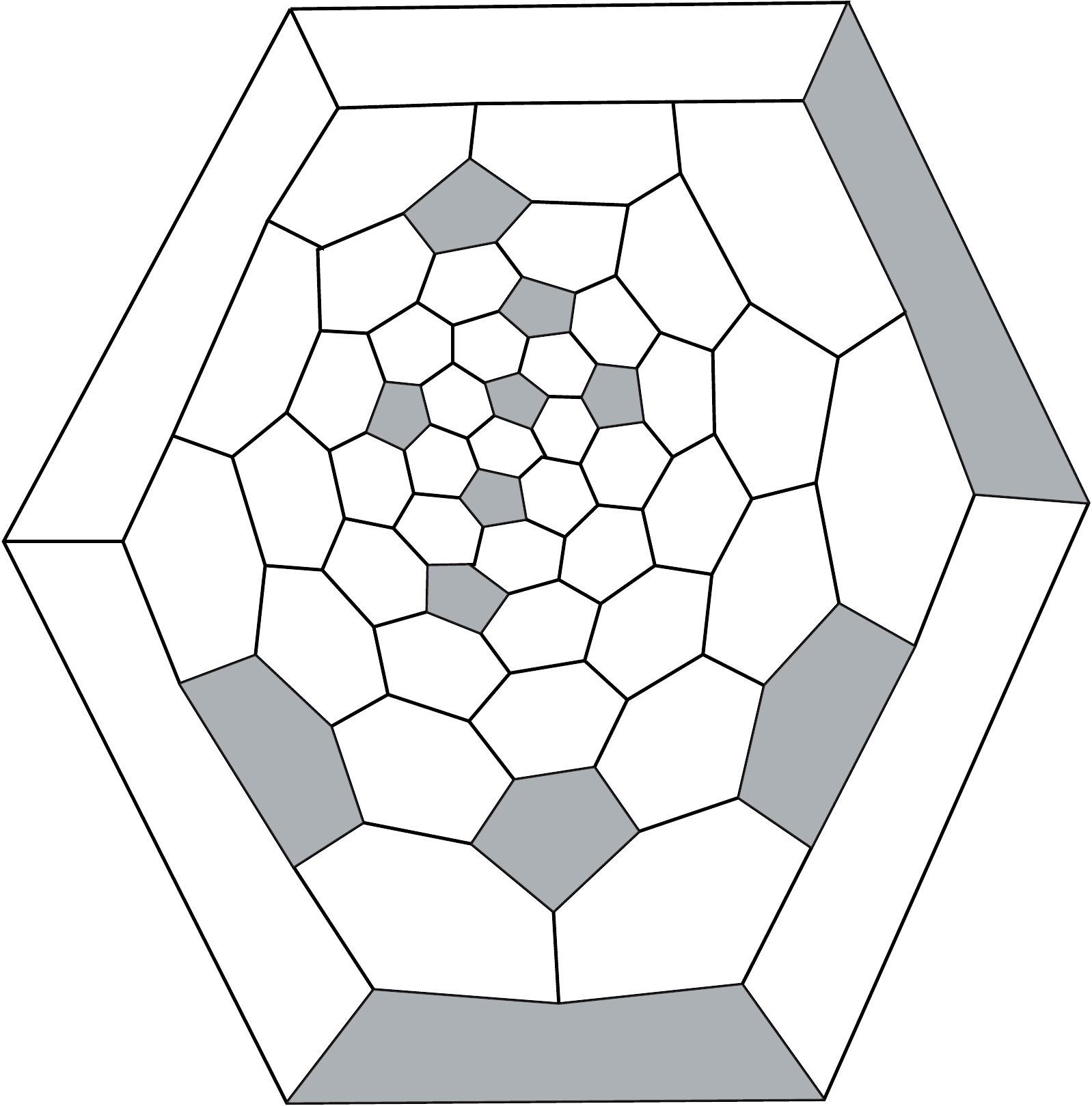}}  $n=110 \ (r=1,2,3, s=3)$ \\
\end{minipage}
\vfill
\begin{minipage}[h]{0.47\linewidth}
\center{\includegraphics[width=0.9\linewidth]{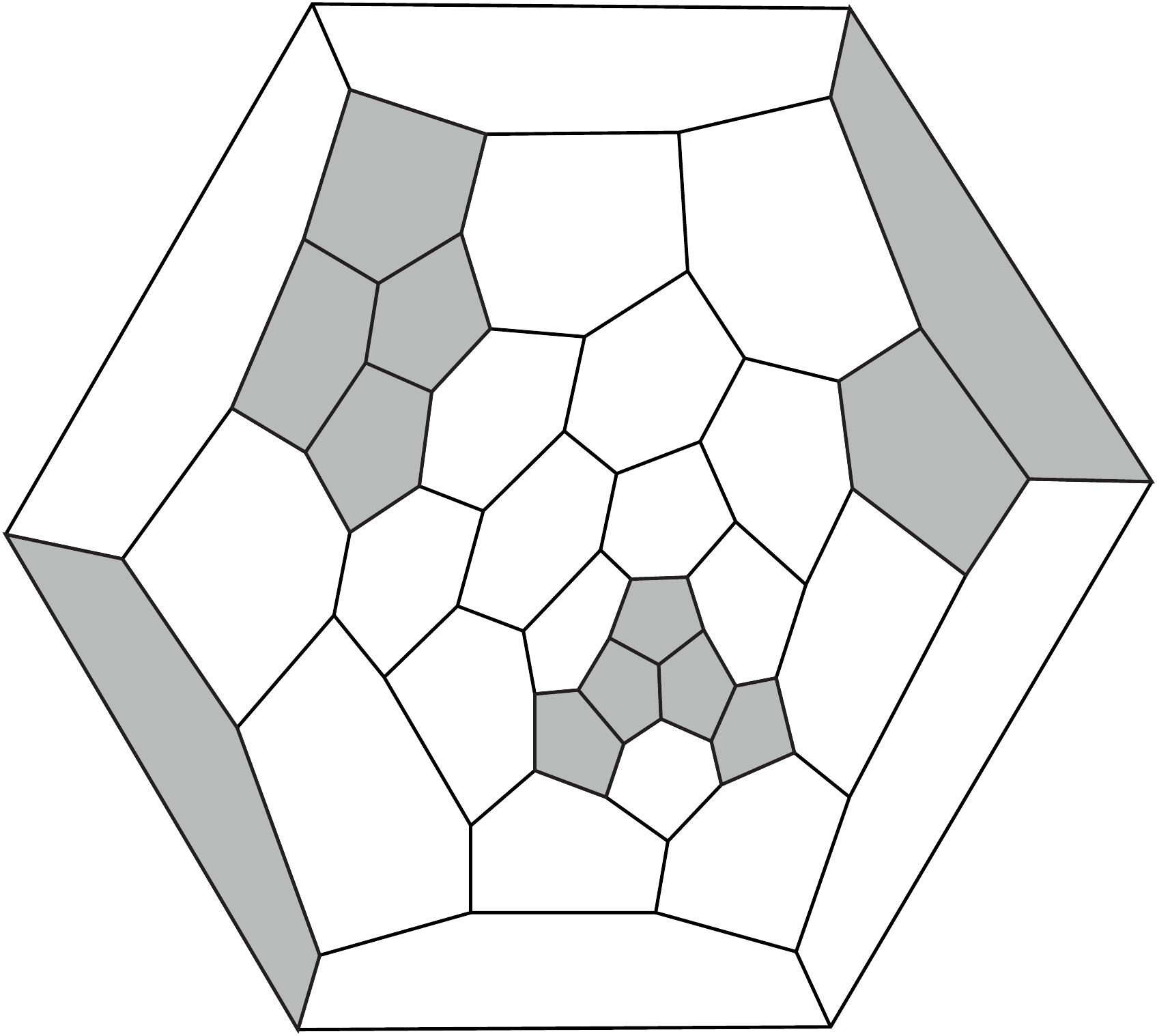}} $n=64 \ (r=s=2)$ \\
\vspace*{5mm}
\end{minipage}
\hfill
\begin{minipage}[h]{0.47\linewidth}
\vspace*{-5mm}
\center{\includegraphics[width=0.85\linewidth]{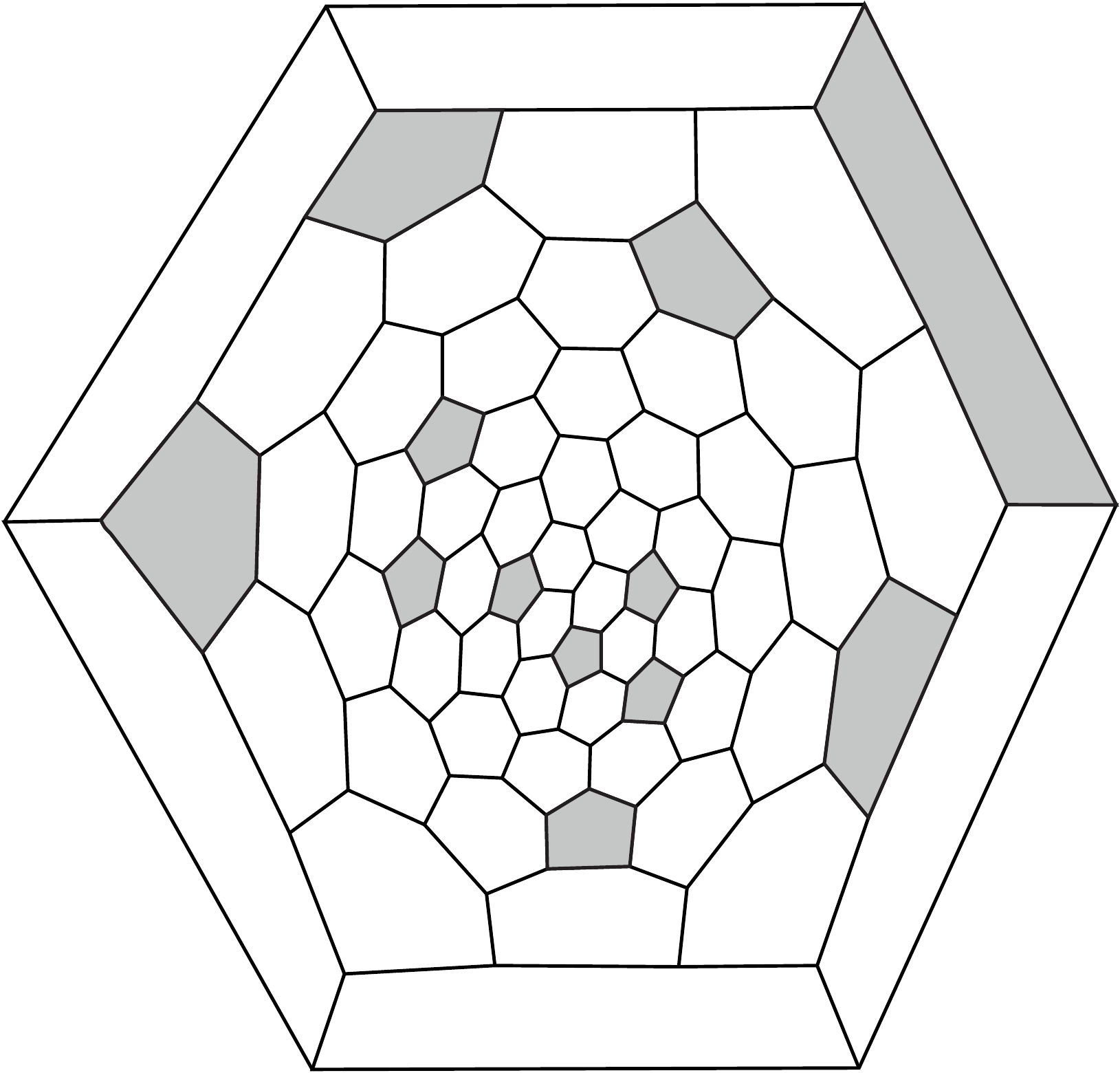}} $n=124 \ (r=s=2)$ \\
\end{minipage}
\vfill
\begin{minipage}[h]{0.47\linewidth}
\center{\includegraphics[width=0.9\linewidth]{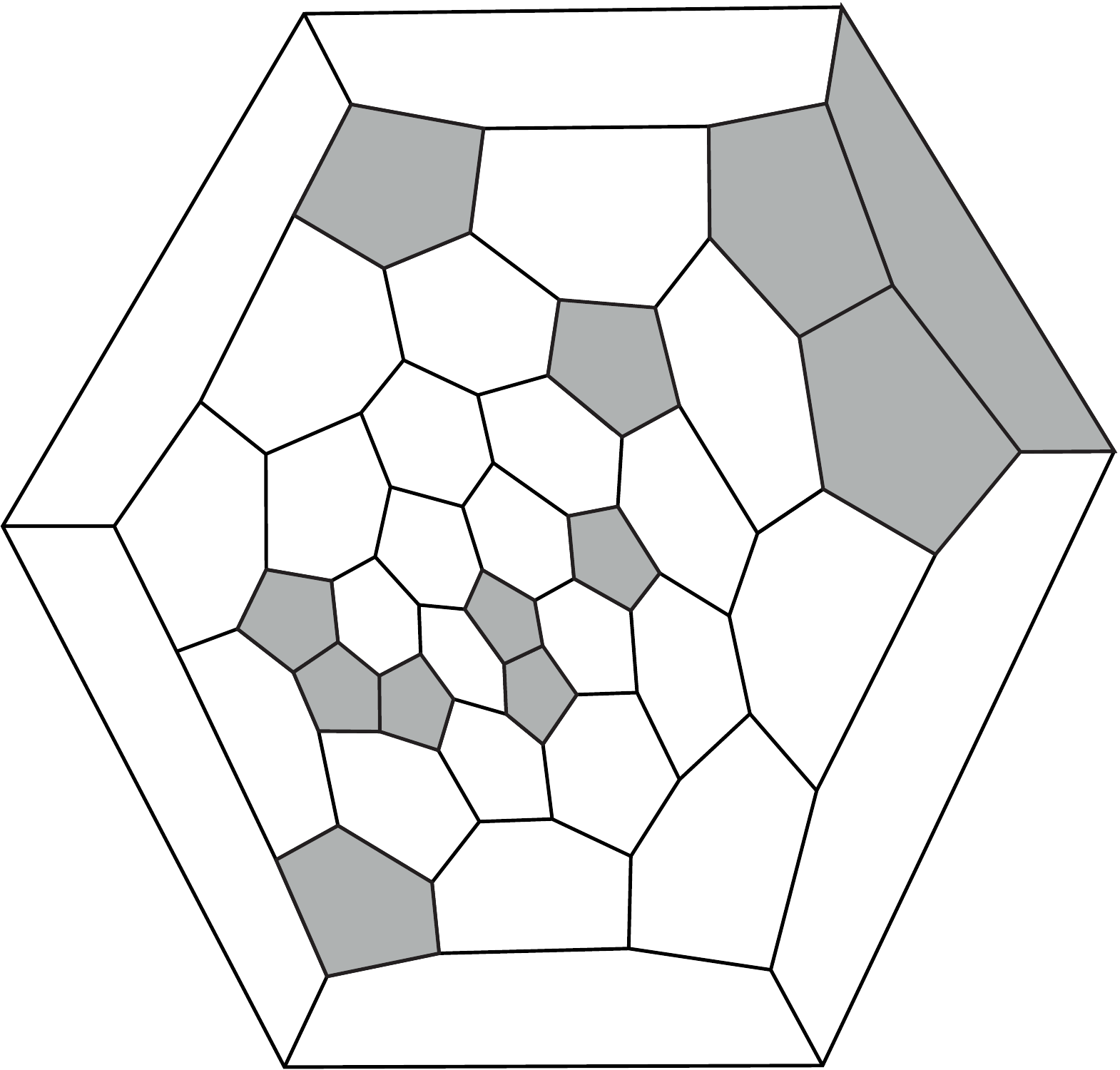}} $n=76 \ (r=1, s=2)$ \\
\vspace*{5mm}
\end{minipage}
\hfill
\begin{minipage}[h]{0.47\linewidth}
\vspace*{-5mm}
\center{\includegraphics[width=0.85\linewidth]{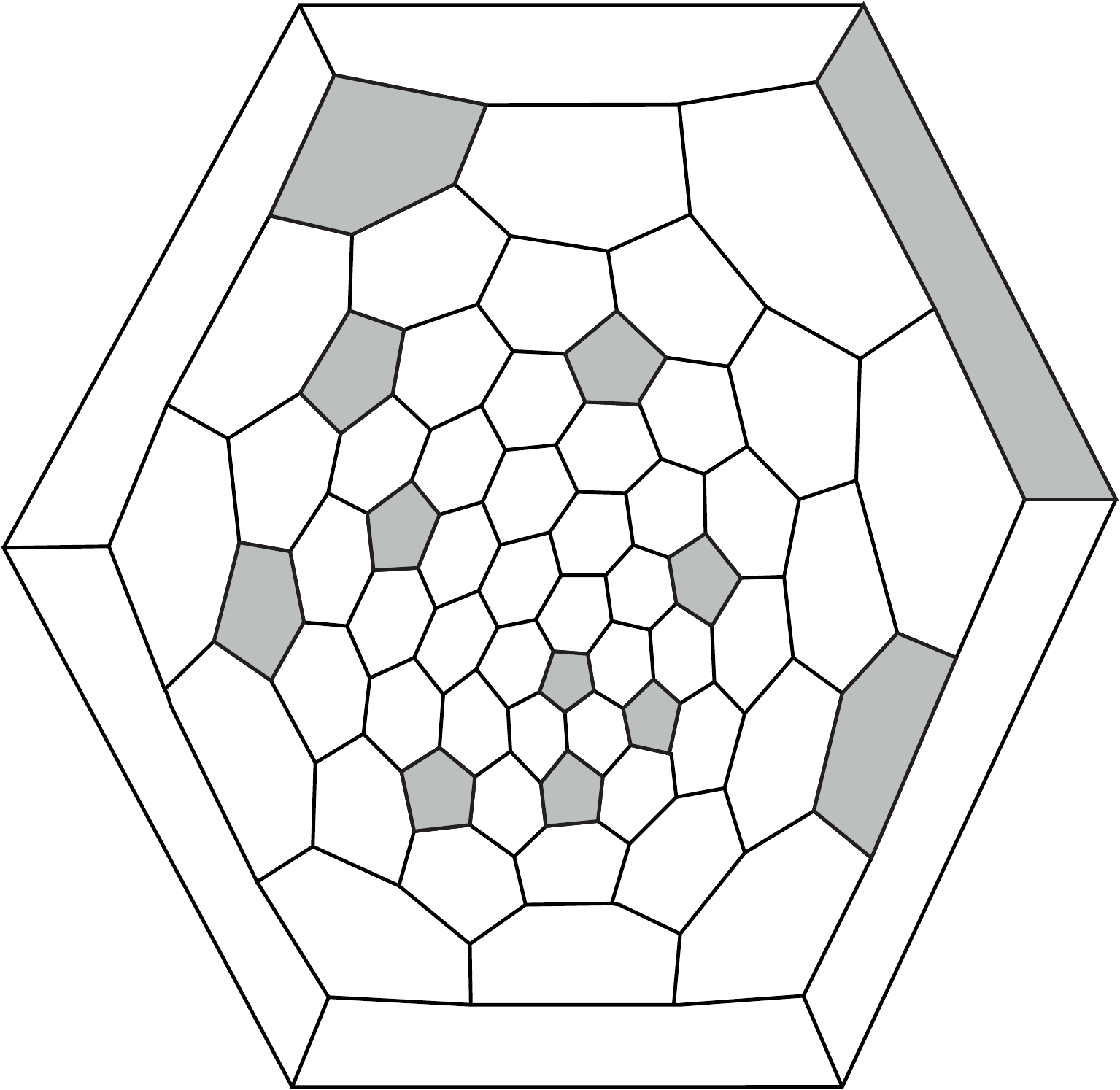}}  $n=130 \ (r=1, s=2)$
\end{minipage}
\caption{Minimal irregular fullerene and IPR fullerene graphs
with $n$ vertices.}
\label{Figure_1}
\end{figure}


Despite the fact that the interval of possible values for transmission
$tr_{1,2}$ is a bit longer than for $tr_{2,2}$,
the number of vertices of irregular fullerene graph for $r=1, s=2$ is less
than for $r=s=2$. The same is also valid for values of $r=s=3$ and $r=2, s=3$.
Based on the growth of irregular fullerene graphs when the number of vertices increases, we
can formulate the following results.

\setcounter{section}{1}

\begin{proposition}
Irregular fullerene graphs with $n$ vertices exist for
$n=64$ if $r=s=2$ or
$$
n \geq \left\{
\begin{array}{rl}
72,& \mbox{if }  $r=s=2$ \\
76,& \mbox{if } $r=1, s=2$ \\
62,& \mbox{if } $r=1,2,3$,\  $s=3$.
\end{array}
\right.
$$
\end{proposition}

\begin{proposition}
Irregular IPR fullerene graphs with $n$ vertices exist for
$$
n \geq \left\{
\begin{array}{rl}
124,& \mbox{if }  $r=s=2$ \\
130,& \mbox{if } $r=1, s=2$ \\
110,& \mbox{if } $r=1,2,3$, \ $s=3$.
\end{array}
\right.
$$
\end{proposition}

Diagrams of the minimal  irregular fullerene graphs are shown in
Fig.~\ref{Figure_1}.
These graphs have 64 vertices ($r=s=2$), 76 vertices ($r=1,s=2$),
and 62 vertices (the same graph for  $r=s=3$,  $r=2, s=3$, and  $r=1, s=3$).
The minimal irregular IPR fullerene graphs
have 124 vertices ($r=s=2$), 130 vertices ($r=1, s=2$),
and 110 vertices (the same graph for  $r=s=3$,  $r=2, s=3$, and  $r=1,s=3$).

As for the Wiener complexity, the following question is still open:

\begin{problem}
Do there exist fullerene graphs with maximum Wiener complexity?
\end{problem}

Based on the above considerations and results of~\cite{Dobr19}, we assume that such irregular fullerene graphs may have 300--400 vertices.

\bigskip

{\bf Acknowledgment.}
The work was supported by the Russian
Foundation for Basic Research (project number 19--01--00682),
the state contract of the Sobolev Institute of Mathematics
(project no. 0314--2019--0016) (AAD),
and the Ministry of Science and Education of Russia
(agreement no. 075--02--2021--1392) (AV).

\end{document}